\documentclass{article}

\usepackage{proof}
\usepackage{latexsym}

\newcommand{\sfj}{{\sf j}}

\newcommand{\blem}{\begin{lemma}}
\newcommand{\elem}{\end{lemma}}
\newcommand{\bth}{\begin{theorem}}
\newcommand{\benu}{\begin{enumerate}}
\newcommand{\eenu}{\end{enumerate}}
\newcommand{\bdes}{\begin{description}}
\newcommand{\edes}{\end{description}}
\newcommand{\bdf}{\begin{definition}}
\newcommand{\edf}{\end{definition}}
\newcommand{\bcor}{\begin{cor}}
\newcommand{\ecor}{\end{cor}}
\newcommand{\bprp}{\begin{proposition}}
\newcommand{\eprp}{\end{proposition}}
\newcommand{\bmlem}{\begin{mlemma}}
\newcommand{\emlem}{\end{mlemma}}
\newcommand{\bclm}{\begin{claim}}
\newcommand{\eclm}{\end{claim}}
\newcommand{\bprf}{{\bf Proof}.\hspace{2mm}}

\newcommand{\eprf}{\hspace*{\fill} $\Box$}

\newcommand{\beqn}{\begin{equation}}
\newcommand{\eeqn}{\end{equation}}
\newcommand{\beqnarr}{\begin{eqnarray}}
\newcommand{\eeqnarr}{\end{eqnarray}}
\newcommand{\beqnarrs}{\begin{eqnarray*}}
\newcommand{\eeqnarrs}{\end{eqnarray*}}

\newcommand{\spand}{\,\&\,}

\newtheorem{theorem}{Theorem}
\newtheorem{definition}[theorem]{Definition}
\newtheorem{proposition}[theorem]{Proposition}
\newtheorem{lemma}[theorem]{Lemma}
\newtheorem{cor}[theorem]{Corollary}
\newtheorem{mlemma}[theorem]{Main Lemma}
\newtheorem{claim}[theorem]{Claim}

\newcommand{\alp}{\alpha}

\newcommand{\veps}{\varepsilon}

\newcommand{\Del}{\Delta}
\newcommand{\ome}{\omega}

\newcommand{\bet}{\beta}
\newcommand{\gam}{\gamma}
\newcommand{\Gam}{\Gamma}

\newcommand{\Sig}{\Sigma}

\newcommand{\Tht}{\Theta}

\newcommand{\vphi}{\varphi}

\newcommand{\fal}{\forall}
\newcommand{\exi}{\exists}

\newcommand{\Rarw }{\Rightarrow}

\newcommand{\lrarw}{\leftrightarrow}
\newcommand{\Lrarw}{\Leftrightarrow}

\newcommand{\setm}{\setminus}

\title{Quick cut-elimination for strictly positive cuts 
}
\author{
Toshiyasu Arai
\thanks{The paper has been finished when I visited M\"unchen. 
I would like to thank to Prof. W. Buchholz for his interests, the valuable comments and the hospitality in my visit.}
\\
Graduate School of Science,
Chiba University
\\
1-33, Yayoi-cho, Inage-ku,
Chiba, 263-8522, JAPAN
\\
tosarai@faculty.chiba-u.jp
}
\date{9.30.2010(M\"unchen)}

\begin{document}
\maketitle
\begin{abstract}
In this paper we show that the intuitionistic theory $\widehat{ID}^{i}_{<\ome}(SP)$ for finitely many iterations of 
strictly positive operators
is a conservative extension of the Heyting arithmetic.
The proof is inspired by the quick cut-elimination due to G. Mints.
This technique is also applied to fragments of Heyting arithmetic.
\end{abstract}

\section{Introduction}\label{sec:intro}

Let us consider in this paper the fixed point predicate $I(x)$ for positive formula $\Phi(X,x)$:
\begin{equation}\label{eq:fix}
(FP)^{\Phi} \; \; \forall x[I(x)  \leftrightarrow  \Phi(I,x)]
\end{equation}

W. Buchholz\cite{Buchholz} showed that an intuitionistic fixed point theory $\widehat{ID}^{i}({\cal M})$ is
conservative over the Heyting arithmetic HA with respect to almost negative formulas(, in which
$\lor$ does not occur and $\exists$ occurs in front of atomic formulas only).
The theory $\widehat{ID}^{i}({\cal M})$ has the axioms (\ref{eq:fix}) $(FP)^{\Phi}$ for fixed points
for {\it monotone formula\/} $\Phi(X,x)$, which is generated from arithmetic atomic formulas and $X(t)$ by means of (first order) monotonic connectives $\lor,\land,\exists,\forall$.
Namely $\to$ nor $\lnot$ does occur in monotone formula.
The proof is based on a recursive realizability interpretation.

After seeing the result of Buchholz, we\cite{attic} showed that
an intuitionistic fixed point (second order) theory is conservative over HA for any arithmetic formulas.
In the theory the operator $\Phi$ for fixed points is generated from  $X(t)$ and 
any second order formulas by means of first order monotonic connectives and second order 
existential quantifiers $\exists f(\in\omega\to\omega)$.
Moreover the same holds for the finite iterations of these operations.
The proof is based on N. Goodman's theorem\cite{Goodman}.

Next C. R\"uede and T. Strahm\cite{Strahm} extends significantly 
the results in \cite{Buchholz} and \cite{attic}.
They showed that the intuitionistic fixed point theory $\widehat{ID}^{i}_{<\ome}(SP)$ for finitely many iterations of 
{\it strictly positive\/} operators
is conservative over HA with respect to negative and $\Pi^{0}_{2}$-formulas.

In this paper we show a full result.

\bth\label{cor:SP}
For each $n<\ome$, $\widehat{ID}^{i}_{n}(SP)$ is conservative over {\rm HA}
with respect to any arithmetic formulas.
In other words
$\widehat{ID}^{i}_{<\ome}(SP)$ is a conservative extension of {\rm HA}.
\end{theorem}

Our proof is based on a quick cut-elimination of strictly positive cuts with arbitrary antecedents,
cf. Theorem \ref{lem:quickelim1}.
The proof is inspired by G. Mints' quick cut-elimination of monotone cuts in  \cite{Mintsmono},
and was found in an attempt to clarify ideas in \cite{IntFx}.

We will give a proof of the non-iterated case, $n=1$, and 
indicate necessary modifications for the general cases
in the subsection \ref{subsec:fiterate}.
Let us explain an idea of our proof more closely. The story is essentially the same as in \cite{IntFx}.
First the finitary derivations in $\widehat{ID}^{i}(SP)$ are embedded to infinitary derivations,
and eliminate cuts partially.
This results in an infinitary derivation of depth less than $\varepsilon_{0}$, and in which there
occurs cut inferences with cut formulas $I^{\Phi}(t)$ for fixed points only.
Now the constraint on operator $\Phi$ admits us to
eliminate strictly positive cut formulas quickly.
In this way we will get an infinitary derivation of depth less than $\varepsilon_{0}$, and in which there
occurs no fixed point formulas.

By formalizing the arguments we see that the end formula is true in HA.

In the section \ref{sec:HM} we show that monotone cuts with negative antecedents 
can be eliminated more quickly.
In the final section \ref{sec:fragments} these techniques are applied to fragments of Heyting arithmetic.

\section{An intuitionistic theory $\widehat{ID}^{i}(SP)$}\label{sec:IDSP}

$L_{HA}$ denote the language of the Heyting arithmetic.
Logical connectives are $\lor,\land,\to, \exists,\forall$.
$\lnot A:\equiv(A\to\bot)$.
Let $I$ be a fresh unary predicate symbol not in $L_{HA}$, and $L_{HA}(I)$ denotes $L_{HA}\cup\{I\}$.

Let $SP$ be the class of $L_{HA}(I)$-formulas
such that $A\in SP$ iff $I$ occurs only strictly positive in $A$.
The class $SP$ is defined inductively.

\begin{definition}
{\rm Define inductively a class of formulas} $SP$ {\rm in} $L_{HA}(I)$
{\rm as follows.}
\begin{enumerate}
\item
{\rm Any atomic formula in} $L_{HA}$ {\rm belongs to} $SP$.
\item
{\rm Any atomic formula} $I(t)$ {\rm belongs to the class} $SP$.
\item
{\rm If} $R,S\in SP${\rm , then} $R\lor S, R\land S, \exists x R,\forall x R\in SP$.
\item\label{df:leftright9}
{\rm If} $L\in L_{HA}$ {\rm and} $R\in SP${\rm , then} $L\to R\in SP$.
\end{enumerate}
\end{definition}

Let $\widehat{ID}^{i}(SP)$ denote the following extension of HA.
Its language is obtained from $L_{HA}$ by adding a unary set constant $I$ for a
$\Phi\equiv\Phi(I,x)\in SP$, in which only a fixed variable $x$ occurs freely.
Its axioms are those of HA in the expanded language(, i.e., the induction axioms are available for
any formulas in the expanded language)  plus the axiom $(FP)^{\Phi}$, (\ref{eq:fix}) for fixed points.

\section{Infinitary derivations}\label{sec:infinitary}

Given an $\widehat{ID}^{i}(SP)$-derivation $D_{0}$ of an $L_{HA}$-sentence $C_{0}$, let us first embed it to
an infinitary derivation in an infinitary calculus $\widehat{ID}^{i\infty}(SP)$.

Let $N$ denote a number which is big enough so that any formula occurring in $D_{0}$ has
logical complexity(, which is defined by the number of occurrences of logical connectives) smaller than $N$.
In what follows any formula occurring in infinitary derivations which we are concerned, has
logical complexity less than $N$.

The derived objects in the calculus $\widehat{ID}^{i\infty}(SP)$ are {\it sequents\/}
$\Gamma\Rightarrow A$, where $A$ is a {\it sentence\/} (in the language of $\widehat{ID}^{i}(SP)$)
and $\Gamma$ denotes a finite set of {\it sentences\/},
where each closed term $t$ is identified with its value $\bar{n}$, the $n$th numeral.

$\bot$ stands ambiguously for false equations $t=s$ with  closed terms $t,s$ having different values.
$\top$ stands ambiguously for true equations $t=s$ with closed terms $t,s$ having same values.

The {\it initial sequents\/} are
\[
\Gamma,I(t)\Rightarrow I(t)\, ; \mbox{\hspace{5mm}} \Gamma,\bot\Rightarrow A\, ; \Gamma\Rightarrow \top
\]

The {\it inference rules\/} are $(L\lor)$, $(R\lor)$, $(L\land)$, $(R\land)$, $(L\to)$, $(R\to)$, 
$(L\exists)$, $(R\exists)$, $(L\forall)$, $(R\forall)$, $(LI)$, $(RI)$ and $(cut)$.
These are standard ones.

\begin{enumerate}

\item
\[
\infer[(LI)]
{\Gamma,I(t) \Rightarrow C}
{\Gamma,\Phi(I,t) \Rightarrow C}
\: ;\:
\infer[(RI)]
{\Gamma\Rightarrow I(t)}
{\Gamma\Rightarrow \Phi(I,t)}
\]
\item
\[
\infer[(L\lor)]
{\Gamma,A_{0}\lor A_{1}\Rightarrow C}
{
\Gamma,A_{0}\Rightarrow C
&
\Gamma,A_{1}\Rightarrow C
}
\: ;\:
\infer[(R\lor)]
{\Gamma\Rightarrow A_{0}\lor A_{1}}
{\Gamma\Rightarrow A_{i}}
\,(i=0,1)
\]
\item
\[
\infer[(L\land)]
{\Gamma,A_{0}\land A_{1}\Rightarrow C}
{
\Gamma,A_{0}\land A_{1},A_{i}\Rightarrow C
}
\, (i=0,1)
\: ;\:
\infer[(R\land)]
{\Gamma\Rightarrow A_{0}\land A_{1}}
{
\Gamma\Rightarrow A_{0}
&
\Gamma\Rightarrow A_{1}
}
\]
\item
\[
\infer[(L\to)]
{\Gamma,A\to B\Rightarrow C}
{
\Gamma,A\to B\Rightarrow A
&
\Gamma,B\Rightarrow C
}
\: ;\:
\infer[(R\to)]
{\Gamma\Rightarrow A\to B}
{\Gamma,A\Rightarrow B}
\]
 \item
\[
\infer[(L\exists)]
{\Gamma,\exists x B(x)\Rightarrow C}
{
\cdots
&
\Gamma,B(\bar{n})\Rightarrow C
&
\cdots (n\in\omega)
}
\: ;\:
\infer[(R\exists)]
{\Gamma\Rightarrow \exists x B(x)}
{\Gamma\Rightarrow B(\bar{n})}
\]
 \item
 \[
 \infer[(L\forall)]
 {\Gamma,\forall x B(x)\Rightarrow C}
 {\Gamma,\forall x B(x),B(\bar{n})\Rightarrow C}
 \: ;\:
 \infer[(R\forall)]
 {\Gamma\Rightarrow \forall x B(x)}
 {
 \cdots
 &
 \Gamma\Rightarrow B(\bar{n})
 &
 \cdots (n\in\omega)
 }
 \]
  \item
  \[
  \infer[(cut)]
  {\Gamma,\Delta\Rightarrow C}
  {
  \Gamma\Rightarrow A
  &
  \Delta,A\Rightarrow C
  }
  \]
\end{enumerate}

The {\it depth\/} of an infinitary derivation is defined to be the depth of the well founded tree.

As usual we see the following proposition.
Recall that $N$ is an upper bound of logical complexities of formulas occurring in the given finite derivation $D_{0}$ of $L_{HA}$-sentence $C_{0}$.

\begin{proposition}\label{prp:embed}
\begin{enumerate}
\item\label{prp:embed1}
There exists an infinitary derivation $D_{1}$ of $C_{0}$ such that its depth is less than $\omega^{2}$ and
the logical complexity of any sentence, in particular cut formulas occurring in $D_{1}$ is less than $N$.
\item\label{prp:embed2}
By a partial cut-elimination, there exist an infinitary derivation $D_{2}$ of $C_{0}$ and
an ordinal $\alpha_{0}<\varepsilon_{0}$
such that
the depth of the derivation $D_{2}$ is less than $\alpha_{0}$ and any cut formula occurring in $D_{2}$ is an atomic formula
$I(t)$(, and the logical complexity of any formula occurring in it is less than $N$).
\end{enumerate}
\end{proposition}

The rank $rk(A)$ of sentences $A$ is defined.

\begin{definition}\label{df:rank1}
{\rm The} rank $rk(A)$ {\rm of a sentence} $A$ {\rm is defined by}
\[
rk(A):= \left\{
\begin{array}{ll}
0 & \mbox{{\rm if} } A\in L_{HA} \\
1 & \mbox{{\rm if }} A\in SP\setm L_{HA} \\
2 & \mbox{{\rm otherwise}}
\end{array}
\right.
\]

\end{definition}

Let us call a cut inference {\it HA-cut\/} [$I${\it -cut\/}] if its cut formula
is of rank 0 [of rank 1], resp.

Let $\vdash^{\alpha}_{r}\Gamma\Rightarrow C$ mean that there exists an infinitary derivation of $\Gamma\Rightarrow C$
such that its depth is at most $\alpha$, and its rank is less than $r$(, and
and the logical complexity of any formula occurring in it is less than $N$).

The following Lemmas are seen as usual.

\begin{lemma}\label{lem:weakfalse}(Weakening lemma)

If $\vdash^{\alp}_{1}\Gamma\Rightarrow A$, then
$\vdash^{\alp}_{1}\Del,\Gam\Rarw A$.
\end{lemma}

\begin{lemma}\label{lem:inversion}(Inversion Lemma)\\
Assume $\vdash^{\alp}_{1}\Gamma\Rightarrow A$.

\begin{enumerate}

\item\label{lem:inversion2}
If $A\equiv B_{0}\land B_{1}$, then $\vdash^{\alp}_{1}\Gamma\Rightarrow B_{i}$ for any $i=0,1$.

\item\label{lem:inversion4}
If $A\equiv \forall x B(x)$, then $\vdash^{\alp}_{1}\Gamma\Rightarrow B(\bar{n})$ for any $n\in\omega$.

\item\label{lem:inversion5}
If $A\equiv B_{0}\to B_{1}$, then $\vdash^{\alp}_{1}\Gamma,B_{0}\Rightarrow B_{1}$.

\item\label{lem:inversion6}
If $A$ is an atomic formula $I(t)$, then $\vdash^{\alp}_{1}\Gamma \Rightarrow \Phi(I,t)$.

\item\label{lem:inversion7}
If $A\equiv\bot$, then
$\vdash^{\alp}_{1}\Gam\Rarw C$ for any $C$.

\item\label{lem:inversion8}
If $\top\in\Gam$, then $\vdash^{\alp}_{1}\Gam_{1}\Rarw A$ for $\Gam_{1}\cup\{\top\}=\Gam$.


\end{enumerate}
\end{lemma}

Let $3_{2}(\bet):=3^{3^{\bet}}$.

\begin{theorem}\label{lem:quickelim1}
Suppose that $\vdash^{\bet}_{2}\Gamma\Rightarrow C$.
Then $\vdash^{3_{2}(\bet)}_{1}\Gamma\Rightarrow C$.
\end{theorem}

Assuming the Theorem \ref{lem:quickelim1}, we can show the Theorem \ref{cor:SP} for $n=1$ as follows.
Suppose an $L_{HA}$-sentence $C_{0}$ is provable in $\widehat{ID}^{i}(SP)$.
By Proposition \ref{prp:embed} we have $\vdash^{\alpha_{0}}_{2}\Rightarrow C_{0}$ for a big enough 
number $N$ and an $\alpha_{0}<\varepsilon_{0}$.
Then Theorem \ref{lem:quickelim1} yields $\vdash^{\beta_{0}}_{1}\Rightarrow C_{0}$
for $\beta_{0}=3_{2}(\alpha_{0})<\varepsilon_{0}$.

Let ${\rm Tr}_{N}(x)$ denote a partial truth definition for formulas of logical complexity less than $N$.
By transfinite induction up to $\beta_{0}$ we see ${\rm Tr}_{N}(C_{0})$.
Note that any sentence occurring  in the witnessed derivation for $\vdash^{\beta_{0}}_{1}\Rightarrow C_{0}$
has logical complexity less than $N$, and it is an $L_{HA}$-sentence.
Specifically there occurs no fixed point formula $I(t)$ in it.
Now since everything up to this point is formalizable in HA, we have ${\rm Tr}_{N}(C_{0})$, and hence $C_{0}$ in HA.
This shows the Theorem \ref{cor:SP} for the case $n=1$.

Additional informations equipped with infinitary derivations 
 together with
 the repetition rule $(Rep)$
  \[
  \infer[(Rep)]
  {\Gamma\Rightarrow C}{\Gamma\Rightarrow C}
  \]
are helpful when we formalize our proof as in \cite{Mintsfinite}.
In this paper let us suppress these.

A proof of Theorem \ref{lem:quickelim1} is given in the next section.

\section{Quick cut-elimination of strictly positive cuts with arbitrary antecedents}\label{sec:proof}

In this section we show that strictly positive cuts can be eliminated quickly
even if antecedents of cut inferences and endsequents are arbitrary formulas.
The only constraint is that any cut formula has to be strictly positive.

Let $\alp\#\bet$ denote the {\it natural sum\/} or commutative sum, $\alp\#\bet=\bet\#\alp$, and
$\alp\times\bet$  the {\it natural product\/}.

Theorem \ref{lem:quickelim1} follows from the following Lemma \ref{lem:main1}.

\blem\label{lem:main1}
For arbitrary $\Gam,\Del$ and $C$,
if
$rk(A)=1$,
\[
\vdash^{\alp}_{1}\Gamma\Rightarrow A
\mbox{ and }
\vdash^{\bet}_{2}\Del,A\Rarw C
\]
then
\[
\vdash^{\alp \times 3_{2}(\bet)}_{1}\Del,\Gam\Rightarrow C
.\]
\elem
{\bf Proof} of Theorem \ref{lem:quickelim1} by induction on $\bet$.
Suppose that $\vdash^{\bet}_{2}\Gamma\Rightarrow C$.
Consider the case when the last rule is an $I$-cut:
\[
\infer[(cut)]{\vdash^{\bet}_{2}\Gamma\Rightarrow C}
{
\vdash^{\gam}_{2}\Gamma\Rightarrow A
&
\vdash^{\gam}_{2}\Gamma,A\Rightarrow C
}
\]
with $rk(A)=1$ and $\gam<\bet$.

By IH(=Induction Hypothesis) we have 
\[
\vdash^{3_{2}(\gam)}_{1}\Gamma\Rightarrow A
\]
Then Lemma \ref{lem:main1} yields 
$\vdash^{3_{2}(\bet)}_{1}\Gamma\Rightarrow C$
since
\[
\gam<\bet\Rarw 3_{2}(\gam)\times 3_{2}(\gam)<3_{2}(\bet)
.\]

This shows Theorem \ref{lem:quickelim1} assuming Lemma \ref{lem:main1}.
\\

\noindent
Next we show Lemma \ref{lem:main1}.
As in Lemma 3.2, \cite{Mintsmono} eliminating procedure is fairly standard,
leaving the resulted cut inferences of rank 0, but has to performed in parallel.

$\mbox{\boldmath$A$}$ denotes a finite list $A_{k},\ldots, A_{2},A_{1}\,(k\geq 0)$ of $SP$-formulas,
and $\mbox{\boldmath$\alp$}=\alp_{k},\ldots,\alp_{2},\alp_{1}$  a list of ordinals.
Then $\vdash^{ \mbox{{\footnotesize \boldmath$\alp$}}}_{1}\Gamma\Rightarrow \mbox{\boldmath$A$}$ 
designates that $\vdash^{\alp_{i}}_{1}\Gam\Rarw A_{i}$ for each $i$.
\[
\sum \mbox{\boldmath$\alp$}:= \left\{
\begin{array}{ll}
\alp_{1}\#\cdots\#\alp_{k} & \mbox{{\rm if }} k>0 \\
1 &  \mbox{{\rm if }} k=0
\end{array}
\right.
\]

$\mbox{\boldmath$A$}_{1}$ denotes the list $A_{k},\ldots, A_{2}$, in which $A_{1}$ is deleted.
Likewise
$\mbox{\boldmath$\alp$}_{1}$ denotes the list $\alp_{k},\ldots, \alp_{2}$.

\blem\label{lem:main1parallel}
Suppose 
$\vdash^{ \mbox{{\footnotesize \boldmath$\alp$}}}_{1}\Gamma\Rightarrow \mbox{\boldmath$A$}$ and 
$\vdash^{\bet}_{2}\Del,\mbox{\boldmath$A$}\Rarw C$.
Then
\beqn\label{eq:main}
\vdash^{(\sum \mbox{\boldmath$\alp$})\times 3_{2}(\bet)}_{1}\Del,\Gam \Rightarrow C
\eeqn
\elem

Note that the case $k=0$ in Lemma \ref{lem:main1parallel} is nothing but Theorem \ref{lem:quickelim1}.

We prove Lemma \ref{lem:main1parallel}
by main induction on $\bet$ with subsidiary induction on $\sum \mbox{\boldmath$\alp$}+k$, where
$k$ is the length of the list $\mbox{\boldmath$A$}$.

\begin{enumerate}
\item\label{1}
The case when one of $\Gam\Rarw A_{i}$, $\Delta,\mbox{\boldmath$A$} \Rightarrow C$ is an initial sequent.

Firs consider the case when $\Delta,\mbox{\boldmath$A$} \Rightarrow C$ is an initial sequent.

If $\Delta,\mbox{\boldmath$A$} \Rightarrow C$ is an initial sequent such that one of the cases $C\equiv\top$, $\bot\in\Delta$ 
or $C\in\Delta$ occurs, then $\Delta\Rightarrow C$, and hence
$\Delta,\Gam\Rightarrow C$
is still the same kind of initial sequent. 

If $\Delta,\mbox{\boldmath$A$} \Rightarrow C$ is an initial sequent with the principal formula $\mbox{\boldmath$A$} \ni A_{i}\equiv C\equiv I(t)$,
then $\Del,\Gamma\Rightarrow A_{i}(\equiv C)$ is an initial sequent.

If $A_{i}\equiv\bot$, then
Inversion lemma \ref{lem:inversion}.\ref{lem:inversion7} with a weakening 
yields $\vdash^{\alp_{i}}_{1}\Del, \Gam\Rarw C$.
$(\sum \mbox{\boldmath$\alp$})\times 3_{2}(\bet)\geq\alp_{i}$ yields (\ref{eq:main}).

Next assume $\Gam\Rarw A_{i}$ is an initial sequent for an $i$.
This implies $k>0$.
For simplicity assume $i=1$.

If $A_{1}\equiv I(t)\in\Gam$, then
we have $\vdash^{ \mbox{{\footnotesize \boldmath$\alp$}}_{1}}_{1}\Gamma\Rightarrow \mbox{\boldmath$A$}_{1}$ and 
$\vdash^{\bet}_{2}\Del,\mbox{\boldmath$A$}_{1},A_{1}\Rarw C$.
If $k=1$, then by weakening $\vdash^{\bet}_{2}\Del,\Gam \Rarw C$, and SIH(=Subsidiary Induction Hypothesis) yields 
$\vdash^{3_{2}(\bet)}_{1}\Del,\Gam \Rarw C$ with $\sum \mbox{\boldmath$\alp$}_{1}=1$ by the definition.

Otherwise by SIH we have 
$\vdash^{(\sum \mbox{\boldmath$\alp$}_{1})\times 3_{2}(\bet)}_{1}\Del,A_{1},\Gam \Rightarrow C$
with $A_{1}\in\Gam$ and $(\sum \mbox{\boldmath$\alp$}_{1})\times 3_{2}(\bet)\leq
(\sum \mbox{\boldmath$\alp$})\times 3_{2}(\bet)$.

If $\bot\in\Gam$, then $\Del,\Gam\Rarw C$ is an initial sequent.

If $A_{1}\equiv\top$, then Inversion lemma \ref{lem:inversion}.\ref{lem:inversion8} 
yields $\vdash^{\bet}_{2}\Del,\mbox{\boldmath$A$}_{1}\Rarw C$, and by SIH
$\vdash^{(\sum \mbox{\boldmath$\alp$}_{1})\times 3_{2}(\bet)}_{1}\Del,\Gam \Rightarrow C$.
\\

\noindent
In what follows assume that none of $\Gam\Rarw A_{i}$, $\Delta,\mbox{\boldmath$A$} \Rightarrow C$ is an initial sequent.

\item\label{3}
Consider the case when $\Delta,\mbox{\boldmath$A$} \Rightarrow C$ is a lowersequent of an $I$-cut.
For a $\gam<\bet$
\[
\infer{\Delta,\Gam\Rightarrow C}
{
\vdash^{\mbox{{\footnotesize \boldmath$\alp$}}}_{1}\Gamma\Rightarrow \mbox{\boldmath$A$} 
&
\infer[(cut)]{\vdash^{\bet}_{2}\Delta,\mbox{\boldmath$A$} \Rightarrow C}
{
 \vdash^{\gam}_{2} \Delta,\mbox{\boldmath$A$} \Rightarrow A_{k+1}
 & 
\vdash^{\gam}_{2} \Delta,\mbox{\boldmath$A$},A_{k+1} \Rightarrow C
 }
}
\]
with $rk(A_{k+1})=1$. 

MIH(=Main Induction Hypothesis) yields $\vdash^{(\sum\mbox{{\footnotesize \boldmath$\alp$}})\times 3_{2}(\gam)}_{1}\Del,\Gam\Rarw A_{k+1}$,
and once again by MIH and 
\[
\left(\sum\mbox{\boldmath$\alp$}\#(\sum\mbox{\boldmath$\alp$})\times3_{2}(\gam)\right) \times 3_{2}(\gam)
\leq
(\sum \mbox{\boldmath$\alp$})\times 3_{2}(\bet)
\]
we conclude $\vdash^{(\sum\mbox{{\footnotesize \boldmath$\alp$}})\times 3_{2}(\bet)}_{1}\Delta,\Gam\Rightarrow C$.

We will depict a `derivation' to illustrate the arguments.
\[
\infer[MIH]{\vdash^{(\sum\mbox{{\footnotesize \boldmath$\alp$}})\times 3_{2}(\bet)}_{1}\Delta,\Gam\Rightarrow C}
{
 \vdash^{\mbox{{\footnotesize \boldmath$\alp$}}}_{1}\Gamma\Rightarrow \mbox{\boldmath$A$} 
 &
\infer{\Del,\Gam,\mbox{\boldmath$A$} \Rarw C}
{
 \infer[MIH]{\vdash^{(\sum\mbox{{\footnotesize \boldmath$\alp$}})\times 3_{2}(\gam)}_{1}\Del,\Gam\Rarw A_{k+1}}
 {
 \vdash^{\mbox{{\footnotesize \boldmath$\alp$}}}_{1}\Gamma\Rightarrow \mbox{\boldmath$A$} 
 &
  \vdash^{\gam}_{2} \Delta,\mbox{\boldmath$A$} \Rightarrow A_{k+1}
  }
 &
 \vdash^{\gam}_{2} \Delta,\mbox{\boldmath$A$},A_{k+1} \Rightarrow C
}
}
\]

In what follows assume that $\Delta,\mbox{\boldmath$A$} \Rightarrow C$ is a lower sequent of an inference rule $J$
other than $I$-cut.

\item\label{4}
If the principal formula of $J$ if any is not in $\mbox{\boldmath$A$} $, then lift up the left upper part:
for a $\gam<\bet$
\[
\infer{\Delta,\Gam\Rightarrow C}
{
\vdash^{\mbox{{\footnotesize \boldmath$\alp$}}}_{1}\Gamma\Rightarrow \mbox{\boldmath$A$} 
&
\infer[(J)]{\vdash^{\bet}_{2}\Delta,\mbox{\boldmath$A$} \Rightarrow C}
{
 \cdots 
 & 
 \vdash^{\gam}_{2} \Delta_{i},\mbox{\boldmath$A$} \Rightarrow C_{i} 
 & 
 \cdots
 }
}
\]

\[
\infer[(J)]
{\vdash^{(\sum \mbox{{\footnotesize \boldmath$\alp$}})\times 3_{2}(\bet)}_{1}\Delta,\Gam\Rightarrow C}
{
\cdots 
&
 \infer[MIH]
 {\vdash^{(\sum\mbox{{\footnotesize \boldmath$\alp$}})\times 3_{2}(\gam)}_{1}\Delta_{i},\Gam\Rightarrow C_{i}}
 {
 \vdash^{\mbox{{\footnotesize \boldmath$\alp$}}}_{1} \Gamma\Rightarrow \mbox{\boldmath$A$}  
 & 
 \vdash^{\gam}_{2}\Delta_{i},\mbox{\boldmath$A$} \Rightarrow C_{i}
 }
&
\cdots
}
\]
Note that $(\sum\mbox{ \boldmath$\alp$})\times 3_{2}(\gam)<(\sum\mbox{ \boldmath$\alp$})\times 3_{2}(\bet)$,
i.e., $(\sum\mbox{ \boldmath$\alp$})>0$ by the definition.

\item\label{5}
Finally suppose that the principal formula of $J$ is a cut formula $A_{i}\in\mbox{\boldmath$A$} $ of $rk(A_{i})=1$.
For simplicity suppose $i=1$.
Use the Inversion Lemma \ref{lem:inversion} if available.
Otherwise examine the left upper part $\vdash^{\alp}_{1}\Gam\Rarw \mbox{\boldmath$A$}$.

 \begin{enumerate}
 \item\label{5a}
The case when $A_{i}\equiv \exists x B(x)\in\mbox{\boldmath$A$} $.

\[
\infer[(L\exists)]
{\vdash^{\bet}_{2}\Delta,\mbox{\boldmath$A$}  \Rightarrow C}
{
 \cdots
 &
 \vdash^{\gam}_{2}\Delta, \mbox{\boldmath$A$} _{1}, B(\bar{n})\Rightarrow C
 &
 \cdots
 }
 \]
 where $A_{1}\not\in\mbox{\boldmath$A$} _{1}$.

We will examine the last rule in $\vdash^{\alp_{1}}_{1}\Gam\Rarw A_{1}(\equiv\exi x B(x))$.
 \benu

 \item\label{5a1}
If $\exi x B(x)$ is derived by an $(R\exi)$,
\[
\infer[(R\exi)]{\vdash^{\alp_{1}}_{1}\Gam\Rarw \exi x B(x)}{\vdash^{\alp_{0}}_{1}\Gam\Rarw B(\bar{n})}
\]
then
 \[
 \infer[MIH]
 {\vdash_{1}^{(\sum\mbox{{\footnotesize \boldmath$\alp$}})\times 3_{2}(\gam)} \Delta,\Gam\Rightarrow C}
 {
   \vdash^{\mbox{{\footnotesize \boldmath$\alp$}}_{1}}_{1}\Gamma\Rightarrow \mbox{\boldmath$A$} _{1}
  &  
  \infer{ \Delta,\Gam, \mbox{\boldmath$A$} _{1}\Rightarrow C}
   {
   \vdash^{\alp_{0}}_{1} \Gamma\Rightarrow B(\bar{n})
   &
   \vdash^{\gam}_{2}\Delta, \mbox{\boldmath$A$} _{1}, B(\bar{n})\Rightarrow C
  }
 }
 \]

\item\label{5a2}
If the last rule is a left rule, then postpone it.


For example $\exi y D(y)\in\Gam$
\[
\infer[(L\exi)]{\vdash^{\alp_{1}}_{1}\Gam\Rarw \exi x B(x)}
{
\cdots
&
\vdash^{\alp_{0}}_{1}\Gam,D(\bar{n})\Rarw \exi x B(x)
&
\cdots
}
\]
Then $\alp_{0}<\alp_{1}$, and hence $\sum\mbox{\boldmath$\alp$}_{1}\#\alp_{0}<\sum\mbox{\boldmath$\alp$}_{1}\#\alp_{1}=\sum\mbox{\boldmath$\alp$}$.
Thus SIH
yields
\[
\vdash^{(\sum \mbox{{\footnotesize \boldmath$\alp$}}_{1}\#\alp_{0})\times 3_{2}(\bet)}_{1}\Del,\Gam,D(\bar{n})\Rarw C\]
 for each $n$.

\[
\hspace{-20mm}
\infer[(L\exi)]{\vdash^{(\sum\mbox{{\footnotesize \boldmath$\alp$}})\times 3_{2}(\bet)}_{1}\Del,\Gam\Rarw C}
{
\cdots
&
 \infer[SIH]{\vdash^{(\sum \mbox{{\footnotesize \boldmath$\alp$}}_{1}\#\alp_{0})\times 3_{2}(\bet)}_{1}\Del,\Gam,D(\bar{n})\Rarw C}
 {
  \vdash^{\mbox{{\footnotesize \boldmath$\alp$}}_{1}}_{1}\Gam\Rarw \mbox{\boldmath$A$}_{1}
 &
  \infer{\Del, \Gam,D(\bar{n}),\mbox{\boldmath$A$}_{1}\Rarw C}
  {
   \vdash^{\alp_{0}}_{1}\Gam, D(\bar{n})\Rarw \exi x B(x)
   &
   \vdash^{\bet}_{2}\Del,\mbox{\boldmath$A$}_{1},\exi x B(x)\Rarw C
   }
 }
&
\cdots
}
\]

Consider next $(L\to)$. Let $D\to E\in\Gam$.
\[
\infer[(L\to)]{\vdash^{\alp_{1}}_{1}\Gam\Rarw \exi x B(x)}
{
\vdash^{\alp_{0}}_{1}\Gam\Rarw D
&
\vdash^{\alp_{0}}_{1}\Gam,E\Rarw \exi x B(x)
}
\]
Then
\[
\hspace{-20mm}
\infer[(L\to)]{\vdash^{(\sum\mbox{{\footnotesize \boldmath$\alp$}})\times 3_{2}(\bet)}_{1}\Del,\Gam\Rarw C}
{
 \vdash^{\alp_{0}}_{1}\Gam\Rarw D
&
 \infer[SIH]{\vdash^{(\sum \mbox{{\footnotesize \boldmath$\alp$}}_{1}\#\alp_{0})\times 3_{2}(\bet)}_{1}\Del,\Gam,E\Rarw C}
 {
  \vdash^{\mbox{{\footnotesize \boldmath$\alp$}}_{1}}_{1}\Gam\Rarw \mbox{\boldmath$A$}_{1}
 &
  \infer{\Del,\Gam,\mbox{\boldmath$A$}_{1}\Rarw C}
  {
  \vdash^{\alp_{0}}_{1}\Gam,E\Rarw \exi x B(x)
  &
  \vdash^{\bet}_{2}\Del, \mbox{\boldmath$A$}_{1}, \exi x B(x) \Rarw C
  }
 }
}
\]

Finally consider an HA-cut with $rk(D)=0$.
\[
\infer[(cut)]{\vdash^{\alp_{1}}_{1}\Gam\Rarw \exi x B(x)}
{
\vdash^{\alp_{0}}_{1}\Gam\Rarw D
&
\vdash^{\alp_{0}}_{1}\Gam,D\Rarw \exi x B(x)
}
\]

\[
\hspace{-20mm}
\infer{\vdash^{(\sum \mbox{{\footnotesize \boldmath$\alp$}})\times 3_{2}(\bet)}_{1}\Del,\Gam\Rarw C}
{
\vdash^{\alp_{0}}_{1}\Gam\Rarw D
&
\infer[SIH]{\vdash^{(\sum \mbox{{\footnotesize \boldmath$\alp$}}_{1}\#\alp_{0})\times 3_{2}(\bet)}_{1}\Del,\Gam,D\Rarw C}
{
     \vdash^{\mbox{{\footnotesize \boldmath$\alp$}}_{1}}_{1}\Gam\Rarw \mbox{\boldmath$A$}_{1}
  &
   \infer{\Del,\Gam,D, \mbox{\boldmath$A$}_{1}\Rarw C}
   {
   \vdash^{\alp_{0}}_{1}\Gam,D\Rarw \exi x B(x)
   &
   \vdash^{\bet}_{2}\Del, \mbox{\boldmath$A$}_{1}, \exi x B(x) \Rarw C
   }
  }
}
\]

\eenu

 \item\label{5b}
 The case when $A_{i}\equiv H\to A_{0}\in\mbox{\boldmath$A$} $ with an $H\in L_{HA}$ and an $A_{0}\in SP$.
For a $\gam<\bet$
\[
\infer[(L\to)]
{\vdash^{\bet}_{2}\Delta,\mbox{\boldmath$A$}  \Rightarrow C}
{
 \vdash^{\gam}_{2}\Delta,\mbox{\boldmath$A$}  \Rightarrow H
&
 \vdash^{\gam}_{2}\Delta,\mbox{\boldmath$A$} _{1},A_{0}\Rightarrow C
}
\]
where $A_{1}\not\in\mbox{\boldmath$A$} _{1}$.

\[
\hspace{-20mm}
\infer[(\mbox{HA-cut})]{\vdash_{1}^{(\sum\mbox{{\footnotesize \boldmath$\alp$}})\times 3_{2}(\bet)}\Delta,\Gam \Rightarrow C}
 {
  \infer[MIH]{\vdash^{(\sum\mbox{{\footnotesize \boldmath$\alp$}})\times 3_{2}(\gam)}_{1}\Delta,\Gam \Rightarrow H}
  {
  \vdash^{\mbox{{\footnotesize \boldmath$\alp$}}}_{1} \Gamma\Rightarrow \mbox{\boldmath$A$} 
   &
  \vdash^{\gam}_{2} \Delta,\mbox{\boldmath$A$} \Rightarrow H
   }
  &
  \infer[MIH]{\vdash^{(\sum\mbox{{\footnotesize \boldmath$\alp$}})\times 3_{2}(\gam)}_{1}\Delta,\Gam,H\Rightarrow C}
  {
  \vdash^{\mbox{{\footnotesize \boldmath$\alp$}}_{1}}_{1}\Gamma\Rightarrow \mbox{\boldmath$A$} _{1}
  &
  \infer{ \Delta,\mbox{\boldmath$A$} _{1},\Gam,H\Rightarrow C}
   {
   \vdash^{\alp_{1}}_{1} \Gamma,H\Rightarrow A_{0}
  &
   \vdash^{\gam}_{2}\Delta,\{A_{0}\}\cup\mbox{\boldmath$A$} _{1} \Rightarrow C
    }
   }
 }
\]
where $\vdash^{\alp_{1}}_{1}\Gamma,H\Rightarrow A_{0}$ by inversion, and $rk(H)=0$.

\item\label{5c}
 The case when $A_{i}\equiv \forall x B(x) \in\mbox{\boldmath$A$} $.
 For a $\gam<\bet$
 \[
\infer[(L\forall)]
{\vdash^{\bet}_{2}\Delta,\mbox{\boldmath$A$}  \Rightarrow C}
{
\vdash^{\gam}_{2}\Delta, \mbox{\boldmath$A$} , B(\bar{n})\Rightarrow C
 }
 \]
By
$(\sum\mbox{\boldmath$\alp$}\#\alp_{1})\times 3_{2}(\gam)\leq(\sum\mbox{ \boldmath$\alp$})\times 3_{2}(\bet)$
and inversion $\vdash^{\alp_{1}}_{1} \Gamma\Rightarrow B(\bar{n})$ we have
 \[
 \infer[MIH]
 {\vdash^{(\sum\mbox{{\footnotesize \boldmath$\alp$}})\times 3_{2}(\bet)}_{1} \Gamma,\Delta\Rightarrow C}
 {
   \vdash^{\mbox{{\footnotesize \boldmath$\alp$}}}_{1}\Gamma\Rightarrow \mbox{\boldmath$A$} 
  &  
  \infer{\Delta,\Gam, \mbox{\boldmath$A$} \Rightarrow C}
   {
   \vdash^{\alp_{1}}_{1} \Gamma\Rightarrow B(\bar{n})
   &
   \vdash^{\gam}_{2}\Delta, \mbox{\boldmath$A$} , B(\bar{n})\Rightarrow C
  }
 }
 \]

\item\label{5d}
 The case when $A_{i}\equiv B_{0}\lor B_{1}\in\mbox{\boldmath$A$} $.
is treated as in the Case (\ref{5a}) for existential quantifier.

\item\label{5e}
 The case when $A_{i}\equiv B_{0}\land B_{1} \in\mbox{\boldmath$A$} $ is treated as in the Case (\ref{5c}) for universal quantifier.
 
\item\label{5f}
The case when $A_{i}$ is a formula $I(t)$.
Use Inversion $\vdash^{\alp}_{1}\Gam\Rarw \Phi(I,t)$.

 \end{enumerate}
\end{enumerate}

This completes a proof of (\ref{eq:main}), and hence of Lemma \ref{lem:main1parallel}.

\subsection{Finite iterations}\label{subsec:fiterate}
Our proof is easily extended to finite iterations of fixed points for strictly positive operators.
The theory $\widehat{ID}^{i}_{n}(SP)$ has the following axiom for formulas $\Phi(X,Y,x)$ in which $X$ occurs only
strictly positive:
\[
\fal i<n\fal x[x\in I^{\Phi}_{i}\lrarw \Phi(I^{\Phi}_{i},I^{\Phi}_{<i},x)]
\]
where $I^{\Phi}_{<i}=\{(y,j): y\in I^{\Phi}_{j}\spand j<i\}$.

Let us explain how to modify the proof.
For simplicity consider the case $n=2$. Drop the superscript $\Phi$ in $I^{\Phi}$, and identify
$I_{<1}$ with $I_{0}$.
Let $\Phi_{i}(X,x):\Lrarw \Phi(X, I_{<i},x)$.
The initial sequents $\Gam,I_{i}(t)\Rarw I_{i}(t)$ and the inference rules $(LI),(RI)$ are for each $i=0,1$
\[
\infer[(LI)_{i}]
{\Gamma,I_{i}(t) \Rightarrow C}
{\Gamma,\Phi_{i}(I_{i},t) \Rightarrow C}
\: ;\:
\infer[(RI)_{i}]
{\Gamma\Rightarrow I_{i}(t)}
{\Gamma\Rightarrow \Phi_{i}(I_{i},t)}
\]

The rank $rk(A)\leq 3$ of sentences $A$ is defined.
Let $SP_{i}$ denote the set of formulas in which $I_{i}$ occurs only strictly positive.
Let $SP_{-1}:=L_{HA}$.

\begin{definition}\label{df:rankn}
{\rm The} rank $rk(A)$ {\rm of a sentence} $A$ {\rm is defined by}
\[
rk(A):= \left\{
\begin{array}{ll}
0 & \mbox{{\rm if }} A\in L_{HA} \\
i+1 & \mbox{{\rm if }} A\in SP_{i}\setm SP_{i-1}\,(i=0,1) \\
3 & \mbox{{\rm otherwise}}
\end{array}
\right.
\]
\edf

Then Theorem \ref{lem:quickelim1} runs as follows.

\begin{theorem}\label{lem:quickelimn}
\benu
\item\label{lem:quickelimn1}
Suppose that $\vdash^{\bet}_{2}\Gamma_{0}\Rightarrow C_{0}$.
Then $\vdash^{3_{2}(\bet)}_{1}\Gamma_{0}\Rightarrow C_{0}$.
\item\label{lem:quickelimn2}
There exists an $m<\ome$ for which the following holds.

Suppose that $\vdash^{\bet}_{3}\Gamma_{0}\Rightarrow C_{0}$.
Then $\vdash^{2_{m}(3_{2}(\bet))}_{2}\Gamma_{0}\Rightarrow C_{0}$.
\eenu
\end{theorem}

Suppose that $\vdash^{\bet}_{3}\Gamma_{0}\Rightarrow C_{0}$.
To prove Theorem \ref{lem:quickelimn}.\ref{lem:quickelimn2}, first eliminate cut inferences with cut formulas
$A$, in which $I_{1}$ does occur strictly positive. 
The proof is the same as in one for Theorem \ref{lem:quickelim1}.
Then the depth of the resulting derivation $D_{1}$ is bounded by $3_{2}(\bet)$.
Unfortunately this derivation $D_{1}$ might be still of rank 2 since for example
if $H\to A_{0}\in SP_{1}$, then $H$ is an arbitrary formula in $L_{HA}(I_{0})$.
In other words $I_{0}$ might occur in $H$ arbitrarily, and we left the cut inference of cut formula $H$
in the Case (\ref{5b}) of the proof of Lemma \ref{lem:main1parallel}.

Now observe that $H$ is a subformula of the fixed operator $\Phi_{1}(I_{1},n)\equiv\Phi(I_{1},I_{0},n)$.
This means that the logical complexity of $H$ can be bounded in advance.
Let $m$ be the number of occurrences  of logical symbols in $\Phi$.
Then in $D_{1}$  eliminate cut inferences of rank 2(, but $I_{1}$ does not occur in their cut formulas),
to get a derivation $D_{2}$ of depth $2_{m}(3_{2}(\bet))$ and of rank 1.
In $D_{2}$ any cut formula is either an HA-formula or of the form $I_{0}(t)$.
Then apply Theorem \ref{lem:quickelimn}.\ref{lem:quickelimn1} to get an HA-derivation of depth
less than $\veps_{0}$.

Obviously this elimination procedure can be iterated finite times as you need, and
the depth of the resulting HA-derivation is less than $\veps_{0}$.
This proves the general case in Theorem \ref{cor:SP}.

\section{Quick cut-elimination of monotone cuts with negative antecedents}\label{sec:HM}

We show that monotone cuts with negative antecedents can be eliminated more quickly.
In this section we consider the Heyting arithmetic HA and its infinitary counterpart $\mbox{{\rm HA}}^{\infty}$.
First let us introduce a class ${\cal NM}$ of $L_{HA}$-formulas.
${\cal N}$ is the class of negative formulas.

\begin{definition}\label{df:leftright}
${\cal N}$ {\rm denotes the class of} negative formulas{\rm , in which no disjunction and existential quantifier occur.}


{\rm Define inductively a class of formulas} ${\cal NM}$ {\rm in} $L_{HA}$
{\rm as follows.}
\begin{enumerate}
\item
{\rm Any atomic formula} $s=t$ {\rm belongs to} ${\cal NM}$.

\item
{\rm If} $R,S\in{\cal NM}${\rm , then} $R\lor S, R\land S, \exists x R,\forall x R\in{\cal NM}$.
\item
{\rm If} $L\in{\cal N}$ {\rm and} $R\in{\cal NM}${\rm , then} $L\to R\in{\cal NM}$.
\end{enumerate}
\end{definition}
It is easy to see ${\cal N}\subset{\cal NM}$.

Note that by the equivalence 
\beqn\label{eq:qntmove}
[\exi x\, A(x)\to B]\lrarw \fal x[A(x)\to B]
\eeqn
$\exi x\, A(x)\to B$ for $A\in{\cal N}$, $B\in{\cal NM}$ is equivalent to the ${\cal NM}$-formula $ \fal x[A(x)\to B]$.

The rank $rk(A)$ of sentences $A$ is redefined as follows.

\begin{definition}\label{df:rank}
{\rm The} rank $rk(A)$ {\rm of a sentence} $A$ {\rm is defined by}
\[
rk(A):= \left\{
\begin{array}{ll}
0 & \mbox{{\rm if} } A\in{\cal N} \\
1 & \mbox{{\rm if }} A\in {\cal NM}\setminus{\cal N} \\
2 & \mbox{{\rm otherwise}}
\end{array}
\right.
\]
\edf

Let $\mbox{HA}^{\infty}$ denote an infinitary system in the language $L_{HA}$, whose initial sequents and inference rules are obtained from those of $\widehat{ID}^{i\infty}(SP)$ by deleting the initial sequents
$\Gamma,I(t)\Rightarrow I(t)$ and inference rules $(LI), (RI)$.

By restricting antecedents to negative (or Harrop) formulas we have a stronger inversion.

\begin{lemma}\label{lem:inversionH}(Inversion Lemma with negative antecedents)\\
Assume $\vdash^{\alp}_{1}\Gamma\Rightarrow A$
such that $\Gamma\subseteq{\cal N}$.

\begin{enumerate}
\item\label{lem:inversionH1}
If $A\equiv B_{0}\lor B_{1}$, then $\vdash^{\alp}_{1}\Gamma\Rightarrow B_{i}$ for an $i=0,1$.

\item\label{lem:inversionH3}
If $A\equiv \exists x B(x)$, then $\vdash^{\alp}_{1}\Gamma\Rightarrow B(\bar{n})$ for an $n\in\omega$.
\end{enumerate}
\end{lemma}

\begin{theorem}\label{lem:quickelimHM}
Let $C_{0}$ denote an ${\cal NM}$-sentence, and $\Gamma_{0}$ a finite set of ${\cal N}$-sentences.
Suppose that $\vdash^{\bet}_{2}\Gamma_{0}\Rightarrow C_{0}$.
Then $\vdash^{2^{\bet}}_{1}\Gamma_{0}\Rightarrow C_{0}$.
\end{theorem}

Again Theorem \ref{lem:quickelimHM} follows from the following Lemma \ref{lem:mainHM} for quick cut-elimination in parallel.

$\mbox{\boldmath$A$}$ denotes a {\it non-empty\/} finite list $A_{k},\ldots, A_{2},A_{1}\,(k> 0)$ of ${\cal NM}$-formulas.
and $\alp$ an ordinal.
Then $\vdash^{\alp}_{1}\Gamma\Rightarrow \mbox{\boldmath$A$}$ 
designates that $\vdash^{\alp}_{1}\Gam\Rarw A_{i}$ for {\it any\/} $i$.
Note here that the depth $\alp$ of the derivations of $\Gam\Rarw A_{i}$ is independent from $i$.

\blem\label{lem:mainHM}
Suppose $\Gam\cup\Del\subset{\cal N}$ and $\mbox{\boldmath$A$}\cup\{C\}\subset{\cal NM}$. 
If
\[
\vdash^{\alp}_{1}\Gamma\Rightarrow \mbox{\boldmath$A$}
\mbox{ and }
\vdash^{\bet}_{2}\Del,\mbox{\boldmath$A$}\Rarw C
\]
then
\[
\vdash^{\alp+2^{\bet}}_{1}\Del,\Gam\Rightarrow C
.\]
\elem

We can prove Lemma \ref{lem:mainHM} by induction on $\bet$ as in Lemma \ref{lem:main1parallel}.
In Case (\ref{1}) we don't need to examine the left upper parts
$\vdash^{\alp}_{1}\Gamma\Rightarrow \mbox{\boldmath$A$}$.
In Case (\ref{5}) the Inversion Lemma on the succedent is always available
since the antecedent $\Gam$ consists solely of negative formulas.
Note that in the Case (\ref{5b}) the remaining cut formula $H\in{\cal N}$ is in the class ${\cal NM}$.

This completes a proof of Lemma \ref{lem:mainHM}, and of Theorem \ref{lem:quickelimHM}.

Note that the procedure leaves cuts with negative cut formulas $H$ in Case (\ref{5b}).
If we restrict to eliminate monotone cuts, then cuts are eliminated quickly and completely.

\begin{theorem}\label{lem:quickelimM}
Let $C_{0}$ denote an ${\cal NM}$-sentence, and $\Gamma_{0}$ a finite set of ${\cal N}$-sentences.
Suppose that there exists a derivation of $\Gamma_{0}\Rightarrow C_{0}$ in which any cut formula is a monotone formula, and whose depth is at most $\bet$.
Then there exists a cut-free derivation of $\Gamma_{0}\Rightarrow C_{0}$ in depth $2^{\bet}$.
\end{theorem}

Let us iterate this procedure for monotone cuts.
\\

\noindent
In what follows $\Phi$ denotes a class of arithmetic formulas such that
any atomic formula is in $\Phi$, and $\Phi$ is closed under substitution of terms for variables and renaming of
bound variables.

Given such a class $\Phi$ of formulas, introduce a hierarchy $\{{\cal M}_{n}(\Phi)\}$ of arithmetic formulas.
\begin{definition}
{\rm First set} ${\cal M}_{1}(\Phi)=\Phi$.

{\rm Define inductively classes of formulas} ${\cal M}_{n+1}(\Phi)\, (n\geq 1)$ {\rm in} $L_{HA}$
{\rm as follows.}
\begin{enumerate}
\item
${\cal M}_{n}(\Phi)\subset{\cal M}_{n+1}(\Phi)$.
\item
{\rm If} $R,S\in {\cal M}_{n+1}(\Phi)${\rm , then} $R\lor S, R\land S, \exists x R,\forall x R\in {\cal M}_{n+1}(\Phi)$.
\item
{\rm If} $L\in {\cal M}_{n}(\Phi)$ {\rm and} $R\in {\cal M}_{n+1}(\Phi)${\rm , then} $L\to R\in {\cal M}_{n+1}(\Phi)$.
\end{enumerate}
\end{definition}
We have $\bigcup_{n<\ome}{\cal M}_{n}(\Phi)=L_{HA}$.

 
For $\Phi=\Sig_{1}$, ${\cal M}_{n}(\Sig_{1})$ coincides with the class $\Tht_{n}$ introduced by W. Burr\cite{Burr}.
Note that by (\ref{eq:qntmove}) for any $n\geq 2$,
each formula in ${\cal M}_{n}(\Sig_{1})=\Tht_{n}$ is equivalent to a formula in  ${\cal M}_{n}(\Del_{0})$,
where $\Del_{0}$ is the class of all atomic formulas.
Also each formula in $\Tht_{2}$ is equivalent to a monotone formula in ${\cal M}$.


The rank $rk(A;\Phi)$ of sentences $A$ relative to the class $\Phi$ is defined.

\begin{definition}\label{df:rank1M}
{\rm The} rank $rk(A;\Phi)$ {\rm of a sentence} $A$ {\rm is defined by}
\[
rk(A;\Phi):=\min\{n-1: A\in{\cal M}_{n}(\Phi)\}
.\]
\end{definition}

Let $\vdash^{\alp}_{r}\Gam\Rarw C$ designate that there exists an infinitary derivation of $\Gam\Rarw C$
such that the depth of the derivation tree is bounded by $\alp$ and any cut formula occurring in it has rank less than 
$r$.
$\vdash^{\alp}_{2}\Gam\Rarw C$ means that in the witnessed derivation of depth $\alp$
any cut formula is in the class ${\cal M}_{2}(\Phi)$.

\begin{theorem}\label{lem:quickelimnM}
Suppose that $\vdash^{\bet}_{r+1}\Gamma_{0}\Rightarrow C_{0}$.
Then $\vdash^{3_{2}(\bet)}_{r}\Gamma_{0}\Rightarrow C_{0}$ for $r\geq 2$.
\end{theorem}
\bprf
This is seen as in the proof of Theorem \ref{lem:quickelimn}.\ref{lem:quickelimn2}, but leave the cut inference of cut formula $H$ with $rk(H;\Phi)<r$
in the Case (\ref{5b}).
\eprf

\section{Applications to fragments of Heyting arithmetic}\label{sec:fragments}

Finally let us remark an application of quick cut-eliminations to fragments of Heyting arithmetic.

\bdf
{\rm Let} $\Phi$ {\rm be a class of arithmetic formulas such that
any atomic formula is in} $\Phi${\rm , and} $\Phi$ {\rm is closed under substitution of terms for variables and renaming of
bound variables.}

$i\Phi$ {\rm denotes the fragment of HA in which induction axioms are restricted to formulas in} $\Phi$.
\[
A(0)\land\fal x[A(x)\to A(x+1)]\to\fal x\, A(x)\, (A\in\Phi)
.\]
{\rm For a class of formulas} $\Psi$,
$\mbox{{\rm RFN}}_{\Psi}(i\Phi)$ {\rm denotes the} $\Psi${\rm -(uniform) reflection principle for} $i\Phi$:
\[
\mbox{{\rm RFN}}_{\Psi}(i\Phi)=\{\mbox{{\rm Pr}}_{i\Phi}(\lceil\vphi(\dot{x})\rceil)\to\vphi(x)
: \vphi\in\Psi\}
\]
{\rm where} $\mbox{{\rm Pr}}_{i\Phi}$ {\rm denotes a standard provability predicate for}
$i\Phi$ {\rm and} $\dot{x}$ {\rm is the} $x${\rm -th formalized numeral.}

{\rm When} $\Psi=L_{HA}$ {\rm the subscript} $\Psi$ {\rm in} $\mbox{{\rm RFN}}_{\Psi}(i\Phi)$ 
{\rm is dropped.}
\edf

By the result of Buchholz\cite{Buchholz} we see that HA proves the consistency of the intuitionistic arithmetic $i{\cal M}$ for the class ${\cal M}$ of monotone formulas 
since $\widehat{ID}^{i}({\cal M})$ can define the truth of monotone formulas, and the consistency statement
$\mbox{{\rm CON}}(i{\cal M})$ is an almost negative formula.
Observe that any prenex $\Pi^{0}_{k}$-formula is a monotone formula, and any monotone formula is equivalent to a prenex formula.

Moreover using truth definition for $\Tht_{n}$-formulas and a partial truth definition
we see
that for each $n\geq 2$
$\widehat{ID}^{i}_{n-1}({\cal M})$ proves the soundness
$\mbox{{\rm RFN}}(i\Tht_{n})$ of $i\Tht_{n}$.
Hence $\mbox{{\rm HA}}\vdash \mbox{{\rm RFN}}(i\Tht_{n})$ by the full conservativity
of $\widehat{ID}^{i}_{n}({\cal M})$ over HA in \cite{attic}.

However this does not show that $\{i{\cal M}_{n}(\Phi)\}_{n}$ forms a proper hierarchy.
Burr\cite{Burr}, Corollary 2.25 shows that $I\Pi^{0}_{n}$ and $i\Tht_{n}$ prove the same $\Pi^{0}_{2}$-sentences
for the fragments $I\Pi^{0}_{n}$ of the Peano arithmetic PA.
Since $I\Pi^{0}_{n+1}$ proves the 2-consistency $\mbox{{\rm RFN}}_{\Pi^{0}_{2}}(I\Pi^{0}_{n})$
of $I\Pi^{0}_{n}$ and hence of $i\Tht_{n}$, by the result of Burr we see that
$i\Tht_{n+1}$ proves the 2-consistency of $i\Tht_{n}$.
Thus $\{i\Tht_{n}\}_{n}$ forms a proper hierarchy.
\\

\noindent
Let us show that $i\Tht_{3}$ proves the soundness of $i\Tht_{2}$ with respect to
$\Tht_{2}$, $\mbox{{\rm RFN}}_{\Tht_{2}}(i\Tht_{2})$.
Recall that $\Tht_{2}$, monotone formulas and formulas in prenex formulas are equivalent each other.

Let $<$ denote a standard $\veps_{0}$-well ordering.
Let
\[
Prg[A] :\Lrarw \fal x[\fal y<x\, A(y) \to A(x)]
\]
and for a class $\Phi$ of formulas, $TI(<\alp,\Phi)$ denote the transfinite induction schema
\[
Prg[A] \to \fal x<\bet \, A(x)
\]
for each $\bet<\alp$ and $A\in\Phi$.

Also let $\ome_{1}:=\ome$ and $\ome_{m+1}:=\ome^{\ome_{m}}$.

\bprp\label{prp:progressive}
If $m+k\leq n+2$, then
\[
i{\cal M}_{n}(\Phi)\vdash TI(<\ome_{m}, {\cal M}_{k}(\Phi))
.\]
\eprp
\bprf
Let
\[
\sfj[A](\alp) :\Lrarw \fal\bet[\fal\gam<\bet A(\gam) \to \fal\gam<\bet+\ome^{\alp} A(\gam)]
.
\]
Then for $A\in{\cal M}_{n}(\Phi)$ we have $\sfj[A]\in{\cal M}_{n+1}(\Phi)$
\[
\mbox{{\rm HA}}({\cal M}_{n}(\Phi))\vdash Prg[A]\to Prg[\sfj[A]]
\]
and $\mbox{{\rm HA}}({\cal M}_{n}(\Phi))\vdash TI(<\ome_{1}, {\cal M}_{n+1}(\Phi))$.
The proposition follows from these.
\eprf

\bcor\label{cor:hierarchyHA}
\benu
\item\label{cor:hierarchyHA1}
For $n\geq 2$
\[
i\Tht_{2n-1}\vdash \mbox{{\rm RFN}}_{\Tht_{2}}(i\Tht_{n})
.\]
For example $i\Tht_{3}$ proves the soundness of prenex induction with prenex consequences.

\item\label{cor:hierarchyHA2} 

For any $m,k,n\geq 1$
\[
i{\cal M}_{2m+k}(\Pi^{0}_{n})\vdash \mbox{{\rm RFN}}_{{\cal M}_{k}(\Pi^{0}_{n})}(i{\cal M}_{m}(\Pi^{0}_{n}))
.\]

\eenu
\ecor
\bprf
\ref{cor:hierarchyHA}.\ref{cor:hierarchyHA1} follows from Theorems \ref{lem:quickelimnM}, \ref{lem:quickelimM}
 and Proposition \ref{prp:progressive}.
Namely embed a finitary derivation of a monotone sentence $C$ in $i{\cal M}_{n}(\Del_{0})$ to an infinitary one.
Apply first Theorem \ref{lem:quickelimnM} $(n-2)$-times, 
to get a derivation of $C$ such that any cut formula occurring in it is a monotone formula and its depth is bounded by
$3_{2n-4}(\ome^{2})=\ome_{2n-3}$.
Then apply Theorem \ref{lem:quickelimM} to get a cut-free derivation of $C$ in depth $2^{\ome_{2n-3}}=\ome_{2n-2}$.
By Proposition \ref{prp:progressive} $TI(<\ome_{2n-1},\Tht_{2})$ is provable in $i\Tht_{2n-1}$.
Since any formula occurring in the cut-free derivation is a subformula of the monotone $C\in\Tht_{2}$,
by a $\Tht_{2}$-truth definition of subformulas of $C$ we knows that $C$ is true in $i\Tht_{2n-1}$.
\\

\noindent
\ref{cor:hierarchyHA}.\ref{cor:hierarchyHA2} follows from Theorem \ref{lem:quickelimnM},
quick cut-elimination of monotone cuts with
arbitrary antecedents and Proposition \ref{prp:progressive}.
Namely embed a finitary derivation of a sentence $C_{0}\in {\cal M}_{k}(\Pi^{0}_{n})$ in 
$i{\cal M}_{m}(\Pi^{0}_{n})$ to an infinitary one.
Eliminate cuts by applying Theorem \ref{lem:quickelimnM} $m$-times, 
and get a derivation of $C_{0}$ in depth $3_{2m}(\ome^{2})=\ome_{2m+1}$,
and in which any cut formula is in $\Pi^{0}_{n}$.
Any formula occurring in the derivation is either a subformula of $C_{0}\in{\cal M}_{k}(\Pi^{0}_{n})$
or a $\Pi^{0}_{n}$-formula.
Therefore using ${\cal M}_{k}(\Pi^{0}_{n})$-truth definition of sequents occurring in the derivation 
and $TI(<\ome_{2m+2},{\cal M}_{k}(\Pi^{0}_{n}))$ we conclude that $C_{0}$ is true in $i{\cal M}_{2m+k}(\Pi^{0}_{n})$.
\eprf
\\

\noindent
Next consider conservations.

The following Corollary \ref{cor:hierarchyHAcons} 
shows, for example that $i\Tht_{2}$ is $\Pi^{0}_{k}$-conservative over $i\Pi^{0}_{k}$
for any $k$, and 
generalizes a theorem by A. Visser and K. Wehmeier(cf. Theorem 3 in \cite{Wehmeier} and Corollary 2.28 in \cite{Burr}.)
stating that
$i\Tht_{2}$ is $\Pi^{0}_{2}$-conservative over $i\Pi^{0}_{2}$.

\bcor\label{cor:hierarchyHAcons}
For any $\Phi\subset\Tht_{2}$,
$i\Tht_{2}$ is $\Phi$-conservative over $i\Phi$.
\ecor
\bprf
Embed a finitary derivation of a monotone sentence $C$ in $i{\cal M}_{2}(\Del_{0})$ to an infinitary one.
Apply Theorem \ref{lem:quickelimM} to get a cut-free derivation of $C$ in depth less than $\ome_{2}$.
\eprf



\end{document}